\documentclass[a4paper]{article}
\usepackage{amsmath}
\usepackage{amsfonts}
\usepackage{amssymb}

\input{tcilatex}
\begin{document}

\begin{center}
{\large On the Gradient of Harmonic Functions\footnote{%
MSC2010:\ Primary 31B05; Secondary 53A07.
\par
Key words:\ harmonic function, gradient, mean curvature, level hypersurface.}%
}

\medskip

Pisheng Ding
\end{center}

\begin{quote}
\textsc{Abstract. }{\small For a harmonic function }$u${\small \ on a domain
in }$%
\mathbb{R}
^{n}${\small , this note shows that }$\left\Vert \nabla u\right\Vert $%
{\small \ is essentially determined by the geometry of level hypersurfaces
of }$u$. {\small Specifically, the factor by which }$\left\Vert \nabla
u\right\Vert ${\small \ changes along a gradient flow\ is completely
determined by the mean curvature of the level hypersurfaces intersecting the
flow.}
\end{quote}

\section{Introduction}

Let $u$ be a harmonic function on an open connected subset $\Omega $ of $%
\mathbb{R}
^{n}$. Suppose, for simplicity, that $u$ has no critical points in $\Omega $%
. Starting at a point $p_{0}\in \Omega $, follow the gradient flow to reach
another point $p\in \Omega $. By how much has $\left\Vert \nabla
u\right\Vert $ changed along the flow? This note seeks to generalize the
answer indicated in \cite{J-R} for the case $n=2$; see Remark 2 in \S 3.

For $a\in u(\Omega )$, let $S_{a}$ denote the level-$a$ hypersurface $%
\left\{ p\mid u(p)=a\right\} $. Orient $S_{a}$ by prescribing the normal
field $\mathbf{N}=\nabla u/\left\Vert \nabla u\right\Vert $. Define $%
H:\Omega \rightarrow 
\mathbb{R}
$ by letting $H(p)$ be the mean curvature of $S_{u(p)}$ at $p$.

\medskip

\noindent \textbf{Theorem}: \textit{Let }$p$\textit{\ be a point on the
gradient flow originating from }$p_{0}$\textit{; let }$\overset{\frown }{%
p_{0}p}$\textit{\ be the arc on the flow between }$p_{0}$\textit{\ and }$p$%
\textit{. Then,}%
\begin{equation*}
\left\Vert \nabla u\left( p\right) \right\Vert =\left\Vert \nabla u\left(
p_{0}\right) \right\Vert \exp \left( (n-1)\int_{\overset{\frown }{p_{0}p}%
}H\,ds\right) \text{ ,}
\end{equation*}%
\textit{where }$s$\textit{\ denotes arc length along }$\overset{\frown }{%
p_{0}p}$\textit{.}

\medskip

It turns out that the case $n=3$ completely embodies the general case. For
simplicity, we treat this case.

\section{Preliminaries}

We begin by considering the mean curvature $H$ of level surfaces of a $C^{2}$
function $f$ on a domain in $%
\mathbb{R}
^{3}$.

Suppose that $0$ is a regular value of $f$ and let $S$ denote the level-$0$
set $\left\{ p\mid f(p)=0\right\} $. Let $\mathbf{N}=\nabla f/\left\Vert
\nabla f\right\Vert $. For $p\in S$, let $T_{p}S$ denote the tangent plane
of $S$ at $p$. For each unit vector $\mathbf{v}\in T_{p}S$, let $\gamma _{%
\mathbf{v}}$ be the unit-speed parametrization of the normal section of $S$
with $\gamma _{\mathbf{v}}(0)=p$ and $\gamma _{\mathbf{v}}^{\prime }(0)=%
\mathbf{v}$, then the (signed) curvature $\kappa _{p}(\mathbf{v})$ of $%
\gamma _{\mathbf{v}}$ at $p$ is defined by the equation $\gamma _{\mathbf{v}%
}^{\prime \prime }(0)=\kappa _{p}(\mathbf{v})\mathbf{N}(p)$. The \textit{%
mean curvature} $H(p)$ of $S$ at $p$ is the mean of $\kappa _{p}(\mathbf{v})$
with $\mathbf{v}$ ranging over the unit circle in $T_{p}S$, i.e.,%
\begin{equation*}
H(p):=\frac{1}{2\pi }\int_{\mathbf{v}\in T_{p}S;\;\left\Vert \mathbf{v}%
\right\Vert =1}\kappa _{p}(\mathbf{v})d\sigma \,\text{,}
\end{equation*}%
where $d\sigma $ is the arc length element of the circle $\left\Vert \mathbf{%
v}\right\Vert =1$.

As $\gamma _{\mathbf{v}}(t)\in S$, the two vectors $\nabla f(\gamma _{%
\mathbf{v}}(t))$ and $\gamma _{\mathbf{v}}^{\prime }(t)$ are orthogonal;
hence%
\begin{eqnarray*}
0 &=&\frac{d}{dt}\left\langle \nabla f(\gamma _{\mathbf{v}}(t)),\;\gamma _{%
\mathbf{v}}^{\prime }(t)\right\rangle \\
&=&\left\langle \frac{d}{dt}\nabla f(\gamma _{\mathbf{v}}(t)),\;\gamma _{%
\mathbf{v}}^{\prime }(t)\right\rangle +\left\langle \nabla f(\gamma _{%
\mathbf{v}}(t)),\;\gamma _{\mathbf{v}}^{\prime \prime }(t)\right\rangle
\end{eqnarray*}%
Note that%
\begin{equation*}
\left\langle \nabla f(\gamma _{\mathbf{v}}(0)),\;\gamma _{\mathbf{v}%
}^{\prime \prime }(0)\right\rangle =\kappa _{p}(\mathbf{v})\left\Vert \nabla
f(p)\right\Vert \text{,}
\end{equation*}%
whereas%
\begin{equation*}
\left\langle \left. \frac{d}{dt}\right\vert _{t=0}\nabla f(\gamma _{\mathbf{v%
}}(t)),\;\gamma _{\mathbf{v}}^{\prime }(0)\right\rangle =Q(p)(\mathbf{v},%
\mathbf{v})\text{,}
\end{equation*}%
where $Q(p):T_{p}%
\mathbb{R}
^{3}\times T_{p}%
\mathbb{R}
^{3}\rightarrow 
\mathbb{R}
$ is the Hessian quadratic form, implying that%
\begin{equation*}
\kappa _{p}(\mathbf{v})=-\frac{1}{\left\Vert \nabla f(p)\right\Vert }Q(p)(%
\mathbf{v},\mathbf{v})\text{.}
\end{equation*}%
To average $\kappa _{p}(\mathbf{v})$ over the unit circle, it suffices to
average $Q(p)(\mathbf{v},\mathbf{v})$. To that end, let $\left( \mathbf{\xi }%
,\mathbf{\eta }\right) $ be an orthonormal basis for $T_{p}S$; then $\left( 
\mathbf{\xi },\mathbf{\eta },\mathbf{N}\right) $ is an orthonormal basis for 
$T_{p}%
\mathbb{R}
^{3}$. It is simple to verify that 
\begin{equation*}
\frac{1}{2\pi }\int_{\mathbf{v}\in T_{p}S\,;\;\left\Vert \mathbf{v}%
\right\Vert =1}Q(p)(\mathbf{v},\mathbf{v})d\sigma =\frac{1}{2}\left(
Q(p)\left( \mathbf{\xi },\mathbf{\xi }\right) +Q(p)\left( \mathbf{\eta },%
\mathbf{\eta }\right) \right) \text{.}
\end{equation*}%
Note that%
\begin{equation*}
Q(p)\left( \mathbf{\xi },\mathbf{\xi }\right) +Q(p)\left( \mathbf{\eta },%
\mathbf{\eta }\right) +Q(p)\left( \mathbf{N},\mathbf{N}\right) =\limfunc{Tr}%
Q(p)=\Delta f\left( p\right) \text{.}
\end{equation*}%
Thus, we have%
\begin{equation}
H(p)=\frac{Q(p)\left( \mathbf{N},\mathbf{N}\right) -\Delta f\left( p\right) 
}{2\left\Vert \nabla f(p)\right\Vert }\text{ .}  \label{Eq Mean Curvature}
\end{equation}

\section{Proof of Theorem}

Let $u$ be a harmonic function on an open connected subset $\Omega $ of $%
\mathbb{R}
^{3}$ without critical points. Let $\mathbf{N}=\nabla u/\left\Vert \nabla
u\right\Vert $. Let $s\mapsto \varphi (s)$ be the unit-speed gradient flow
originating from $p_{0}$; $\varphi $ is such that%
\begin{equation*}
\varphi ^{\prime }(s)=\mathbf{N}(\varphi (s))\text{\quad and\quad }\varphi
(0)=p_{0}\text{.}
\end{equation*}

Consider $g(t):=u(\varphi (t))$. Then,%
\begin{eqnarray*}
g^{\prime }(s) &=&\left\langle \nabla u(\varphi (s)),\;\varphi ^{\prime
}(s)\right\rangle =\left\Vert \nabla u(\varphi (s))\right\Vert \text{,} \\
g^{\prime \prime }(s) &=&\left\langle \frac{d}{ds}\nabla u(\varphi
(s),\;\varphi ^{\prime }(s)\right\rangle +\left\langle \nabla u(\varphi
(s)),\;\varphi ^{\prime \prime }(s)\right\rangle \text{.}
\end{eqnarray*}%
As $\left\Vert \varphi ^{\prime }\right\Vert \equiv 1$, the vectors $\varphi
^{\prime }$ and $\varphi ^{\prime \prime }$, and hence $\nabla u$ and $%
\varphi ^{\prime \prime }$, are always orthogonal. Thus,%
\begin{equation*}
g^{\prime \prime }(s)=Q(\varphi (s))(\mathbf{N},\mathbf{N})\text{.}
\end{equation*}%
By (\ref{Eq Mean Curvature}) in \S 2,%
\begin{eqnarray}
H(\varphi (s)) &=&\frac{Q(\varphi (s))\left( \mathbf{N},\mathbf{N}\right) }{%
2\left\Vert \nabla u(\varphi (s))\right\Vert }=\frac{1}{2}\frac{g^{\prime
\prime }(s)}{g^{\prime }(s)}=\frac{1}{2}\frac{d}{ds}\log g^{\prime }(s) 
\notag \\
&=&\frac{1}{2}\frac{d}{ds}\log \left\Vert \nabla u(\varphi (s))\right\Vert 
\text{.}  \label{Eq Directional Derivative}
\end{eqnarray}%
Integrating both sides yields the result.

\bigskip

\noindent \textbf{Remarks}:

\begin{enumerate}
\item If $u$ is a harmonic function of two variables, the above analysis can
be easily adapted to show that%
\begin{equation*}
\left\Vert \nabla u\left( p\right) \right\Vert =\left\Vert \nabla u\left(
p_{0}\right) \right\Vert \exp \left( \int_{\overset{\frown }{p_{0}p}}\kappa
\,ds\right) \text{,}
\end{equation*}%
where $\kappa $ is the curvature (signed, according to choice made of $%
\mathbf{N}$) of the level curves of $u$. This formula was suggested in \cite%
{J-R} by a purely complex-analytic argument, as $u$ is locally the real part
of a holomorphic function. However, this argument does not apply when $n>2$.

\item For $n>3$, it suffices to notes that a level hypersurface of an $n$%
-variable $C^{2}$ function $f$ has at $p$ mean curvature%
\begin{equation*}
H(p)=\frac{Q(p)\left( \mathbf{N},\mathbf{N}\right) -\Delta f\left( p\right) 
}{(n-1)\left\Vert \nabla f(p)\right\Vert }\text{ ,}
\end{equation*}%
where the factor $1/(n-1)$, the counterpart of the factor $1/2$ in (\ref{Eq
Mean Curvature}), stems from averaging a quadratic form on $%
\mathbb{R}
^{n}$ over a unit $(n-2)$-sphere. Omitting the obvious details, we conclude,
as a generalization of (\ref{Eq Directional Derivative}), that, for a
harmonic function $u$ of $n$ variables,%
\begin{equation*}
H=\frac{1}{(n-1)}D_{\mathbf{N}}\log \left\Vert \nabla u\right\Vert \text{ ,}
\end{equation*}%
from which the general case of the Theorem follows.
\end{enumerate}

Illinois State University, Normal, Illinois

\texttt{pding@ilstu.edu}

\end{document}